\documentclass[A4paper,twoside,11pt,english]{article}
\usepackage{babel}
\usepackage{graphicx}
\usepackage{latexsym}
\usepackage{amsmath}
\usepackage{amsthm}
\usepackage{amsfonts}
\usepackage{mathrsfs}
\usepackage{geometry}
\usepackage[latin1]{inputenc}
\usepackage{amssymb}
\usepackage{makeidx}
\usepackage{xcolor}
\usepackage{makeidx}
\usepackage{tocbibind}
\usepackage{tcolorbox}
\usepackage[framemethod=TikZ]{mdframed}
\usepackage{lipsum}
\usepackage[mathscr]{eucal}
\usepackage{array}
\usepackage{cite}
\usepackage{hyperref}

\title{\textbf{Advancements in nonlinear exponential sampling: convergence, quantitative analysis and Voronovskaya-type formula}}
\author{\textbf{Danilo Costarelli}\thanks{Corresponding author} \quad - \quad \textbf{Mariarosaria Natale} \vspace{0.3cm}\\
Department of Mathematics and Computer Science \\
            University of Perugia\\
        1, Via Vanvitelli, 06123 Perugia, Italy    \\  
{\small {\tt danilo.costarelli@unipg.it}} \quad - \quad  {\small {\tt mariarosaria.natale@unipg.it}}}
\date{}

\providecommand{\U}[1]{\protect\rule{.1in}{.1in}}
\theoremstyle{plain}
\newtheorem{theorem}{Theorem}

\newtheorem{corollary}[theorem]{Corollary}

\newtheorem{definition}[theorem]{Definition}

\newtheorem{lemma}[theorem]{Lemma}

\theoremstyle{definition}
\newtheorem{remark}[theorem]{Remark}

\def\R{{\mathbb{R}}}
\def\Z{{\mathbb{Z}}}
\def\supp{{\text{supp}}}

\begin{document}
\maketitle
\begin{abstract}
{\small \noindent In this paper, we introduce the nonlinear exponential Kantorovich sampling series. We establish pointwise and uniform convergence properties and a nonlinear asymptotic formula of the Voronovskaja-type given in terms of the limsup. Furthermore, we extend these convergence results to Mellin-Orlicz spaces with respect to the logarithmic (Haar) measure. Quantitative results are also given, using the log-modulus of continuity and the log-modulus of smoothness, respectively, for log-uniformly continuous functions and for functions in Mellin-Orlicz spaces. Consequently, the qualitative order of convergence can be obtained in case of functions belonging to suitable Lipschitz (log-H\"olderian) classes.}
\end{abstract}

\section{Introduction}
The analysis of nonlinear operators is a very current and dynamic field of study, bridging pure mathematics and practical applications in technology and science. A notable example is in Signal Processing, where nonlinear transformations, during signal filtering, produce new frequencies. This nonlinear behaviour is also evident in power electronics and wireless communications, which involve amplifiers that introduce nonlinear distortions to their input signals. Similarly, in radiometric photography and CCD image sensors, the relationship between input radiance and intensity exhibits nonlinearity, even though it is monotone increasing. Amplifier saturation also introduces nonlinear distortions into the input signal. The pioneering works in the theory of nonlinear operators, particularly in connection with approximation problems, and a wide literature can be found in \cite{Musielak1983bis, Musielak1991, Musielak1993, MonografiaBMV, VZ2022, CCNV2022, CNV2023, CNVnew}.\\

Concurrently, significant progress has been made in the study of approximation results by means of exponential-type operators (see, e.g., \cite{ACDarielli2024, acar2023pointwise, aldaz2009bernstein, angeloni2024approximation, acar2024bivariate, aral2022generalized, angamuthu2020direct, bardaro2021durrmeyer}). In this respect, Bardaro, Faina, and Mantellini \cite{bardaro2017faina} introduced a new family of sampling-type operators, known as exponential sampling series. Their approach uses exponentially spaced sample values to represent Mellin band-limited functions, proving particularly useful in optical physics phenomena, such as light scattering and Fraunhofer diffraction \cite{bertero1991exponential, casasent2006optical, gori1993sampling, ostrowsky1981exponential}. From a mathematical perspective, these operators are best studied within the framework of Mellin analysis, particularly using Mellin transform theory (see \cite{butzer1998exponential, butzer1997direct, boccuto2022some}).\\

Motivated by these developments, Costarelli in 	\cite{chapter_costa}, proposed a generalization that extends to the nonlinear framework the operators introduced by Bardaro, Faina, and Mantellini in 2017. This generalization is expressed as
\begin{equation*}
(S^\chi_w f)(x):=\sum_{k\in\mathbb{Z}}\chi\left ( e^{-t_{k}}x^w, f(e^{t_k/w}) \right ), \;\;\;\; x\in\mathbb{R}^{+},
\end{equation*}
where $f:\R^{+}\to\R$ is any function such that the above series is convergent for every $x\in\R^+$, and where the bivariate function $\chi:\R^{+}_{0}\times\R\to\R$ is a nonlinear kernel and satisfies suitable assumptions.\\

In this paper, we further this study by presenting an analogous generalization for the Kantorovich exponential sampling series, based on the linear version described in \cite{BBSV2007}, that here takes the following form
\begin{equation*}
(K^\chi_w f)(x):=\sum_{k\in\mathbb{Z}}\chi\left ( e^{-t_{k}}x^w,\frac{w}{\Delta_{k}}\int_{\frac{t_k}{w}}^{\frac{t_{k+1}}{w}}f(e^u)du \right ), \;\;\;\; x\in\mathbb{R}^{+},
\end{equation*}
where the sample values of the function $f$ are replaced with their mean values over small intervals (see Definition \ref{DefOP}). Such operators are suitable for proving approximation results in Mellin Lebesgue spaces, or more in general, in Mellin-Orlicz or modular spaces, with the aim to obtain approximation results in the case of not-necessarily log-uniformly continuous signals.\\

Here we obtain pointwise and uniform convergence theorems, also giving a quantitative version in terms of the so-called log-modulus of continuity (see Section \ref{MainSec}), and a nonlinear asymptotic formula of the Voronovskaja type under certain local regularity assumptions on the function $f$. Similarly, the qualitative order of approximation has been established in the case of functions belonging to suitable locally log-H\"{o}lderian classes. The last section (Section \ref{SectMellO}) is devoted to extend the convergence properties within the more general framework of Mellin-Orlicz spaces via a density approach. First we prove a modular convergence theorem for the nonlinear exponential Kantorovich sampling operators acting on the space of continuous functions with compact support, then we obtain a modular-type inequality, and finally we exploit a well-known density result for log-uniformly continuous function with compact support in modular spaces. Finally, we estimate the rate of convergence (Section 5) through both quantitative and qualitative analysis in the setting of Mellin-Orlicz spaces.

%---------------PRELIMINARIES--------------------------------------------
\section{Basic notions and preliminaries}\label{SECpreliminary}
%Let us denote by $\mathbb{N}$ the set of positive integers and by $\Z$ the set of all integers. Similarly, we denote by $\R$ and $\R^{+}$the sets of all real and positive real numbers, respectively. 
Throughout this paper, $C^0(\R^+)$ will represent the space of all continuous and bounded functions defined on $\R^{+}$. A function $f\in C^0(\R^+)$ is said to be \textit{log-uniformly continuous} on $\R^+$, if for any positive $\varepsilon$, there exists a $\delta>0$ such that for all $u,v\in\R^+$,
$$|f(u)-f(v)|<\varepsilon \text{ whenever } |	\ln u -\ln v|\leq \delta.$$ 
We denote by $\mathcal{C}(\R^+)$ the subspace of $C^0(\R^+)$ consisting of log-uniformly continuous (and bounded) functions. It is important to note that a log-uniformly continuous function is not necessarily uniformly continuous, and vice versa. However, these two notions are equivalent on compact intervals within $\R^+$.

Furthermore, the class of all Lebesgue measurable functions and the space of all bounded function on $\R^+$ will be denoted by $M(\R^+)$ and $L^{\infty}(\R^+)$, respectively. Finally, let $\mathcal{M}(\R^+)$ and $L^{1}_{\mu}(\R^{+})$ be respectively the set of all Haar measurable functions and the Lebesgue space with respect to the logarithmic (Haar) measure $\mu$. We recall that, $\mu$ is defined for any set $A\subset \R^{+}$ as
\begin{equation*}
	\mu(A):=\int_{A}\frac{dt}{t}.
\end{equation*}

Now, we recall the definition of the class of operators we work with. \\
Let $\Pi=(t_{k})_{k\in\Z}$ be a sequence of real numbers with $-\infty<t_{k}<t_{k+1}<+\infty$, $\lim_{k\to \pm \infty}t_{k}=\pm \infty$, and such that there exist $ \Delta, \delta >0$ for which $\delta\leq \Delta_{k}:=t_{k+1}-t_{k}\leq \Delta$, for every $k\in\Z$. 

A function $\chi:\R^+\times\R\to\R$ will be called a \textit{nonlinear kernel} (see \cite{chapter_costa}) if it satisfies the following conditions
\begin{description}
\item ($\chi$1) $\left( \chi(e^{-t_k}x^w,u)\right)_{k}\in \ell^{1}(\Z)$, for every $x\in\R^+$, $u\in\mathbb{R}$ and $w>0$;
\item ($\chi$2) $\chi(x,0)=0$, for every $x\in\R^+$;
\item ($\chi$3) $\chi$ is an $(L,\psi)$-Lipschitz kernel, i.e., there exist a measurable function $L:\R^+\to\R^{+}_{0}$ and a non-decreasing and continuous function $\psi:\R^+_0\to\R^{+}_0$, with $\psi(0)=0$, such that 
\begin{equation*}
\left | \chi(x,u)-\chi(x,v) \right |\leq L(x)\psi (\left | u-v \right |),	
\end{equation*}
for every $x\in\R^{+}$ and $u,v\in\R$;
\item ($\chi$4) there exists $\alpha>0$ such that, for every $j\in\mathbb{N}$ and $w>0$,
\begin{enumerate}
	\item \begin{equation*}
	\mathcal{S}_{w}^{j}(x):=\sup_{0\leq |u|<\frac{1}{j}}\left | \sum _{k\in\mathbb{Z}}\chi \left (e^{-t_k}x^w,u \right )-u \right |= \mathcal{O}\left (w^{-\alpha}  \right ),
	\end{equation*}
\item \begin{equation*}
\mathcal{T}^{j}_w(x):=\sup_{\frac{1}{j}\leq |u|} \left |\frac{1}{u} \sum _{k\in\mathbb{Z}}\chi \left ( e^{-t_k}x^w,u \right )-1 \right |= \mathcal{O}\left (w^{-\alpha}  \right ),
\end{equation*}
\end{enumerate}
as $w\to +\infty$, uniformly with respect to $x\in\R^{+}$.\newline
\end{description}
\noindent Moreover, we assume that the function $L$ of condition $(\chi 3)$ satisfies the following additional assumptions
\begin{description}
\item ($L$1) $L\in L^1_{\mu}(\mathbb{R}^+)$ and is bounded;
\item ($L$2) there exists a number $\beta>0$ such that
\begin{equation*}
M_{\beta , \Pi}(L):=\sup_{x\in\mathbb{R}^{+}}\sum_{k\in\mathbb{Z}}L(e^{-t_k}x^w) \left |  \ln x-t_{k}\right |^{\beta}<+\infty ,	
\end{equation*}
i.e., the log-discrete absolute moment of order $\beta$ is finite.
\end{description}

\begin{definition}\label{DefOP}
We define the nonlinear exponential sampling operators of the Kantorovich-type for any given nonlinear kernel $\chi$ as follows
\begin{equation*}
(K^\chi_w f)(x):=\sum_{k\in\mathbb{Z}}\chi\left ( e^{-t_{k}}x^w,\frac{w}{\Delta_{k}}\int_{\frac{t_k}{w}}^{\frac{t_{k+1}}{w}}f(e^u)du \right ), \quad x\in\mathbb{R}^{+},\quad w>0,
\end{equation*}
for any given locally integrable function $f:\mathbb{R}^{+}\to\mathbb{R}$ such that the above series is convergent, for every $x\in\mathbb{R}^{+}$.
\end{definition}
These operators represent a Kantorovich version of the nonlinear exponential sampling series, introduced in \cite{chapter_costa}. They are derived by replacing sample values of the form $f(e^{t_{k}/w})$ with Steklov integrals of order 1 of the form $\frac{w}{\Delta_{k}}\int_{\frac{t_k}{w}}^{\frac{t_{k+1}}{w}}f(e^u)du$ (see, e.g., \cite{costa_new_steklov}), for locally integrable function $f$ defined on $\R^{+}$. Through the change of variables $u=\ln v +t_{k}/w$, $K_{w}^{\chi}f$ can be expressed as follows
\begin{equation*}
(K^\chi_w f)(x)=\sum_{k\in\mathbb{Z}}\chi\left ( e^{-t_{k}}x^w,\frac{w}{\Delta_{k}}\int_{1}^{e^{\Delta_k/w}}f(ve^{t_k/w})\frac{dv}{v} \right ).
\end{equation*}
\begin{remark}
Note that, if $\chi(x, u) := L(x) u$, where $L$ satisfies the conditions ($L$1) and ($L$2), the operators $K_w^\chi f$ reduce to the linear exponential  Kantorovich sampling operators considered in \cite{acar2023pointwise, acar2024bivariate, aral2022generalized, angamuthu2020direct}.
\end{remark}
We now recall the following lemma that will be useful in the next sections.
\begin{lemma}[Lemma 1.1 of \cite{chapter_costa}]\label{lemma_preliminare}
Under the above assumptions on $L$, we have:
\begin{description}
\item[(i)] $M_{0,\Pi}(L) := \displaystyle\sup_{x \in \mathbb{R}_+}  \sum_{k \in \mathbb{Z}} L(e^{-t_k}x) < +\infty $;
\item[(ii)]  for every $\gamma > 0$, it turns out that
    \begin{equation*}
        \sum_{|t_k-w\ln x|>\gamma w} L(e^{-t_k}x^w) \leq \frac{1}{\gamma^\beta w^{\beta}} M_{\beta,\Pi}(L),
    \end{equation*}
    for every $x\in\R^+$, from which we get:
    \begin{equation*}
        \lim_{w \to +\infty}  \sum_{|t_k-w\ln x|>\gamma w} L(e^{-t_k}x^w) = 0,
    \end{equation*}
    uniformly with respect to $x \in \mathbb{R}^+$.
\end{description}
\end{lemma}

\begin{remark}\label{remark_bounded}
If we assume that $f\in L^{\infty}(\R^{+})$, by condition ($\chi$3) and by (i) of Lemma \ref{lemma_preliminare}, we obtain that $(K^{\chi}_w f)_{w>0}$ are well-defined. In fact, it turns out that
\begin{equation*}
	\begin{split}
		|(K^\chi_w f)(x)|&\leq\sum_{k\in\mathbb{Z}}\left|\chi\left ( e^{-t_{k}}x^w,\frac{w}{\Delta_{k}}\int_{\frac{t_k}{w}}^{\frac{t_{k+1}}{w}}f(e^u)du \right )\right|\\
		&\leq \sum_{k\in\mathbb{Z}} L(e^{-t_{k}}x^w) \psi\left(\frac{w}{\Delta_{k}}\int_{\frac{t_k}{w}}^{\frac{t_{k+1}}{w}}|f(e^u)|du\right)\\
		&\leq M_{0,\Pi}(L) \psi(\left\|f\right\|_{\infty})<+\infty,
	\end{split}
\end{equation*}
for every $x\in\mathbb{R}^{+}$ and $w>0$.
\end{remark}

\begin{remark}\label{remark_no_f_bounded}
The assumption that $f$ is bounded on $\R^{+}$ in Remark \ref{remark_bounded} can be relaxed. Instead of requiring strict boundedness, it suffices to assume the existence of two positive constants $a$ and $b$, such that $|f(e^x)|\leq a+b|x|$, for every $x\in\R^{+}$, i.e., $f$ behaves like a logarithm as the argument approaches infinity provided that certain absolute moments are finite.
\\
Let $x,w\in\R^{+}$ be fixed. Suppose that condition ($\chi 3$) is satisfied with $L$ such that $M_{1,\Pi}(L)<+\infty$ and $\psi(u)=u$, $u\in\R^{+}$ (i.e., the case in which a strongly Lipschitz condition is satisfied).  By the growth condition on $f$, we have
\begin{equation*}
\begin{split}
	|(K^\chi_w f)(x)|&\leq \sum_{k\in\mathbb{Z}} L(e^{-t_{k}}x^w) \frac{w}{\Delta_{k}}\int_{\frac{t_k}{w}}^{\frac{t_{k+1}}{w}}	\left(a+b|u|\right)du\\
	&\leq (a+b|\ln x|)M_{0,\Pi}(L)+b\sum_{k\in\mathbb{Z}} L(e^{-t_{k}}x^w) \frac{w}{\Delta_{k}}\int_{\frac{t_k}{w}}^{\frac{t_{k+1}}{w}}	|u-\ln x|du\\
	&\leq \left(a+b|\ln x|+\frac{b\Delta}{2w}\right)M_{0,\Pi}(L)+\frac{b}{w}M_{1,\Pi}(L)<+\infty.
\end{split}
\end{equation*}
\end{remark}
%--------------MAIN SECTION ---------------------------------------------
\section{Results in $\mathcal{C}(\R^{+})$}\label{MainSec}
We can prove the following convergence theorem.
\begin{theorem}\label{PointUnifTH}
	Let $f:\R^+\to\R$ be a bounded function. Then,
\begin{equation*}
	\lim_{w\to +\infty} (K^\chi_w f)(x)=f(x).
\end{equation*}
at any point $x\in\mathbb{R}^+$ of continuity of $f$. Moreover, if $f\in \mathcal{C}(\mathbb{R}^+)$, then
\begin{equation*}
	\lim_{w\to +\infty}\left\|K^\chi_w f - f\right\|_{\infty}=0.
\end{equation*}
\end{theorem}
\begin{proof}
We prove only the first part of the theorem, since the second one can be obtained by similar arguments. Let $x\in\R^+$ be a fixed point of continuity of $f$. We estimate the error of approximation $|(K^\chi_w f)(x)- f(x)|$, obtaining by ($\chi$3)
\begin{align}
|(K^\chi_w f)(x)- f(x)|&\leq \sum_{k\in\Z}\left|\chi\left ( e^{-t_{k}}x^w,\frac{w}{\Delta_{k}}\int_{\frac{t_k}{w}}^{\frac{t_{k+1}}{w}}f(e^u)du \right )-\chi\left ( e^{-t_k}x^w,f(x)\right )\right|\notag\\
&+\left|\sum_{k\in\Z}\chi\left ( e^{-t_k}x^w,f(x)\right )-f(x)\right|\notag\\
&\leq \sum_{k\in\Z} L( e^{-t_k}x^w)\psi\left(\frac{w}{\Delta_{k}}\int_{\frac{t_k}{w}}^{\frac{t_{k+1}}{w}}|f(e^u)-f(x)|du\right)\notag\\
&+\left|\sum_{k\in\Z}\chi\left ( e^{-t_k}x^w,f(x)\right )-f(x)\right|\notag\\
&=:I_1+I_2.\label{split_cont}
\end{align}
We estimate $I_1$. By the continuity of $f$ at $x$, for every fixed $\varepsilon>0$ there exists $\gamma>0$ such that $|f(x)-f(e^u)|=|f(e^{\ln x})-f(e^u)|<\varepsilon$ whenever $\left|\ln x -u\right|\leq\gamma$. Now, we split $I_1$ in two additional summands, namely $I_1=I_{1,1}+I_{1,2}$, where
\begin{equation*}
	I_{1,1}:= \sum_{\left|t_k-w\ln x\right|\leq \frac{\gamma w}{2}} L( e^{-t_k}x^w)\psi\left(\frac{w}{\Delta_{k}}\int_{\frac{t_k}{w}}^{\frac{t_{k+1}}{w}}|f(e^u)-f(x)|du\right),
\end{equation*}
and 
\begin{equation*}
	I_{1,2}:= \sum_{\left|t_k-w\ln x\right|> \frac{\gamma w}{2}} L( e^{-t_k}x^w)\psi\left(\frac{w}{\Delta_{k}}\int_{\frac{t_k}{w}}^{\frac{t_{k+1}}{w}}|f(e^u)-f(x)|du\right).
\end{equation*}
For every $u\in \left[\frac{t_k}{w},\frac{t_{k+1}}{w}\right]\subset \R$, if $\left|t_k-w\ln x\right|\leq \frac{\gamma w}{2}$, we have
\begin{equation}\label{disGamma}
	\left|u-\ln x\right|\leq \left|u-\frac{t_{k}}{w}\right|+\left|\frac{t_{k}}{w}-\ln x\right|\leq \frac{\Delta}{w}+\frac{\gamma}{2}\leq \gamma,
\end{equation}
since we can choose $w>0$ sufficiently large to satisfy $\frac{\Delta}{w}\leq\frac{\gamma}{2}$. Hence,
\begin{equation*}
I_{1,1}\leq \sum_{\left|t_k-w\ln x\right|\leq \frac{\gamma w}{2}} L(e^{-t_k}x^w)\psi\left(\frac{w}{\Delta_{k}}\int_{\frac{t_k}{w}}^{\frac{t_{k+1}}{w}}\varepsilon \,du\right)\leq M_{0,\Pi}(L)\,\psi(\varepsilon),
\end{equation*}
for sufficiently large $w>0$. For $I_{1,2}$, there holds
\begin{equation*}
	I_{1,2}\leq \psi(2\left\|f\right\|_{\infty}) \sum_{\left|t_k-w\ln x\right|> \frac{\gamma w}{2}} L(e^{-t_k}x^w).
\end{equation*}
By (ii) of Lemma \ref{lemma_preliminare}, it follows that $I_{1,2}\to 0$ as $w\to+\infty$, uniformly with respect to $x\in\R^+$. Since $\varepsilon>0$ is arbitrary and $\psi$ is continuous, then also $I_1\to 0$ as $w\to+\infty$.\newline
Now, we estimate $I_2$. Let $A_j:=\{x\in\R^+\,:\,0\leq |f(x)|<1/j\}$, with $j\in\mathbb{N}$ fixed,  taking into account condition ($\chi$2), we can rewrite $I_2$ as follows
\begin{equation*}
\begin{split}
I_2&\leq\left|\sum_{k\in\Z}\chi\left (e^{-t_k}x^w,f(x)\mathbf{1}_{A_j}\right )-f(x)\mathbf{1}_{A_j}\right|\\
&+\left|\sum_{k\in\Z}\chi\left (e^{-t_k}x^w,f(x)\mathbf{1}_{\R^+\setminus A_j}\right )-f(x)\mathbf{1}_{\R^+\setminus A_j}\right|\\
&\leq \left|\sum_{k\in\Z}\chi\left (e^{-t_k}x^w,f(x)\mathbf{1}_{A_j}\right )-f(x)\mathbf{1}_{A_j}\right|+|f(x)|\mathcal{T}^{j}_w(x)\\
&=:I_{2,1}+I_{2,2}.
\end{split}
\end{equation*}
It is clear that $I_{2,1}\leq \mathcal{S}_{w}^{j}(x)$ and $I_{2,2}\leq \left\|f\right\|_{\infty}\mathcal{T}^{j}_w(x)$, hence using ($\chi$4) we have that $I_{2}\to 0$ as $w\to+\infty$ uniformly with respect to $x\in\R^+$. This concludes the proof.
\end{proof}
To establish quantitative estimates for the nonlinear exponential Kantorovich sampling operators, we begin by recalling the notion of the \textit{log-modulus of continuity} for a function $f : \mathbb{R}^+ \rightarrow \mathbb{R}$. This is defined as:
\begin{equation*}
\widetilde{\omega}(f,\delta) := \sup \left\{ |f(x) - f(y)| : |\ln x - \ln y| \leq \delta, \, x, y \in \mathbb{R}^+ \right\}, 
\end{equation*}
for any positive parameter $\delta > 0$. This notion was first introduced in \cite{bardaro2014mellin}, and it satisfies all the properties of the classical modulus of continuity (see, \cite{bardaro2017faina}). Specifically, if $f\in \mathcal{C}(\R^+)$, then $\widetilde\omega(f,\delta)\to 0$ as $\delta\to 0^+$.\\

It is important to note that the estimate presented in Theorem \ref{stima_c} below depends on whether condition (L2) holds with $\beta \geq 1$ or $0 < \beta < 1$. Indeed, there are kernels, such as the Mellin-Fejer kernel (see, e.g., \cite{DO2017, CH2018, chapter_costa}), for which the discrete absolute moments of order \(\beta \geq 1\) are not finite, while condition (L2) is satisfied for some \(0 < \beta < 1\). Therefore, we provide the following theorem by addressing these two cases separately.

\begin{theorem}\label{stima_c}
Let $\chi$ be a nonlinear kernel satisfying ($\chi$3) with $\psi$ concave, and let $f\in \mathcal{C}(\R^+)$. Then, we have the following inequalities, for sufficiently large $w>0$
\begin{enumerate}
	\item if $L$ satisfies condition (L2) with $\beta\geq 1$:
\begin{equation*}
\left\|K^\chi_w f-f\right\|_{\infty}\leq M_3\psi\left(\widetilde{\omega}\left(f,\frac{1}{w}\right)\right)+ M_1w^{-\alpha}+M_2\left\|f\right\|_{\infty} w^{-\alpha},
\end{equation*}
\item if $L$ satisfies condition ($L$2) with $0<\beta< 1$:
\begin{equation*}
\begin{split}
\left\|K^\chi_w f-f\right\|_{\infty}&\leq M_4\,\psi\left(\widetilde\omega\left(f, w^{-\beta}\right)\right)+ 2^{\beta+1}\;\psi(\left\|f\right\|_{\infty})\;w^{-\beta}\;M_{\beta,\Pi}(L)\\
&+ M_1 w^{-\alpha}+M_2\left\|f\right\|_{\infty}  w^{-\alpha},
\end{split}
\end{equation*}
\end{enumerate}
where $M_3:=M_{0,\Pi}(L)+\Delta\, M_{0,\Pi}(L)+M_{1,\Pi}(L)$, $M_4:=M_{0,\Pi}(L)+M_{\beta,\Pi}(L)+\,\Delta^{\beta}M_{0,\Pi}(L)$, and $M_1,M_2,\alpha>0$ are the constants of condition ($\chi 4$).
\end{theorem}
\begin{proof} 
Let $x \in \mathbb{R}^+$ be fixed. The approximation error $|(K^\chi_w f)(x) - f(x)|$ can be decomposed by $I_1+I_2$, as illustrated in the preliminary steps of the proof of Theorem \ref{PointUnifTH} (see (\ref{split_cont})).\\
Now, we focus on estimating $I_1$. We have to distinguish between two cases: when the discrete absolute moment of order $\beta$ is finite for $\beta\geq 1$ (Case 1), and when it is finite for $0<\beta<1$ (Case 2).
\begin{description}
	\item{\textit{(Case 1).}}  Since $\psi$ is non decreasing, we get 
\begin{equation*}
	\begin{split}
		I_1&\leq \sum_{k\in\Z} L(e^{-t_k}x^w)\,\psi\left(\frac{w}{\Delta_{k}}\int_{\frac{t_k}{w}}^{\frac{t_{k+1}}{w}}\widetilde{\omega}(f,\left|u-\ln x\right|)du\right)\\
		&\leq \sum_{k\in\Z} L(e^{-t_k}x^w)\,\psi\left(\frac{w}{\Delta_{k}}\int_{\frac{t_k}{w}}^{\frac{t_{k+1}}{w}}\widetilde{\omega}\left(f,\frac{1}{w}\right)\left[1+w\left|u-\ln x\right|\right]du\right)\\
		&= \sum_{k\in\Z} L(e^{-t_k}x^w)\,\psi\left(\widetilde{\omega}\left(f,\frac{1}{w}\right)\left[1+\frac{w}{\Delta_{k}}\int_{\frac{t_k}{w}}^{\frac{t_{k+1}}{w}}w\left|u-\ln x\right|du\right]\right),
	\end{split}
\end{equation*}
for every $w>0$, where the previous estimate is a consequence of $\widetilde{\omega}(f,\lambda\delta)\leq (1+\lambda)\widetilde{\omega}(f,\delta)$, with $\lambda=w\left|u-\ln x\right|$ and $\delta=\frac{1}{w}$. Now, for every $x\in\R^+$, and $u\in\R$, we may write
\begin{equation}\label{diseq}
	\begin{split}
		\left|u-\ln x\right|\leq \left|u-\frac{t_{k}}{w}\right|+\left|\frac{t_{k}}{w}-\ln x\right|\leq \frac{\Delta}{w}+\frac{\left|w\ln x-t_{k}\right|}{w},
	\end{split}
\end{equation}
for every $w>0$; therefore 
\begin{equation*}
	I_1\leq \sum_{k\in\Z} L(e^{-t_k}x^w)\,\psi\left(\widetilde{\omega}\left(f,\frac{1}{w}\right)\left[1+\Delta+\left|w\ln x-t_{k}\right|\right]\right).
\end{equation*}
Since $\psi$ is concave, we have for $u\geq 1$
\begin{equation}\label{disPSI}
	u\psi(v)=u\psi\left(\frac{1}{u}\cdot vu \right)\geq u\frac{1}{u}\psi(vu)=\psi(vu),
\end{equation}
for every $v\geq 0$; consequently, we finally get
\begin{equation*}
\begin{split}
		I_1&\leq \sum_{k\in\Z} L(e^{-t_k}x^w)\,\left[1+\Delta+\left|w\ln x-t_{k}\right|\right]\psi\left(\widetilde{\omega}\left(f,\frac{1}{w}\right)\right)\\
		&\leq M_{0,\Pi}(L)\;(1+\Delta)\;\psi\left(\widetilde{\omega}\left(f,\frac{1}{w}\right)\right)+\sum_{k\in\Z} L(e^{-t_k}x^w)\,\left|w\ln x-t_{k}\right|\psi\left(\widetilde{\omega}\left(f,\frac{1}{w}\right)\right)\\
		&\leq M_{0,\Pi}(L)\;(1+\Delta)\;\psi\left(\widetilde{\omega}\left(f,\frac{1}{w}\right)\right)+M_{1,\Pi}(L)\;\psi\left(\widetilde{\omega}\left(f,\frac{1}{w}\right)\right).
\end{split}
\end{equation*}
\item{\textit{(Case 2).}} We now split the series in $I_1$ in the following manner
\begin{equation*}
\begin{split}
I_1&\leq \biggl\{\sum_{|w\ln x-t_k|\leq w/2}+\sum_{|w\ln x-t_k|> w/2}\biggr\} L(e^{-t_k}x^w)\,\psi\left(\frac{w}{\Delta_{k}}\int_{\frac{t_k}{w}}^{\frac{t_{k+1}}{w}}|f(e^u)-f(x)|du\right)\\
&=:I_{1,1}+I_{1,2}.
\end{split}
\end{equation*}
Before estimating $I_{1,1}$, we observe that, for every $u\in \left[\frac{t_k}{w},\frac{t_{k+1}}{w}\right]$, and if $|w\ln x-t_k|\leq w/2$ we have
\begin{equation*}
		\left|u-\ln x\right|\leq \left|u-\frac{t_{k}}{w}\right|+\left|\frac{t_{k}}{w}-\ln x\right|\leq \frac{\Delta}{w}+\frac{1}{2}\leq 1,
\end{equation*}
for $w>0$ sufficiently large, and moreover, since $0<\beta<1$, it is also easy to see that
\begin{equation*}
\widetilde\omega(f, |u-\ln x|)\leq \widetilde\omega(f,|u-\ln x|^{\beta}).
\end{equation*}
Hence, by using again the property for which $\widetilde\omega(f,\lambda \delta)\leq (1+\lambda)\widetilde\omega(f,\delta)$, with $\lambda=(w\left|u-\ln x\right|)^{\beta}$ and $\delta=w^{-\beta}$, we get 
\begin{equation*}
\begin{split}
I_{1,1}&\leq \sum_{|w\ln x-t_k|\leq w/2}L(e^{-t_k}x^w)\,\psi\left(\frac{w}{\Delta_{k}}\int_{\frac{t_k}{w}}^{\frac{t_{k+1}}{w}}\widetilde\omega\left(f, |u-\ln x|^{\beta}\right)du\right)\\
&\leq  \sum_{|w\ln x-t_k|\leq w/2}L(e^{-t_k}x^w)\,\psi\left(\frac{w}{\Delta_{k}}\int_{\frac{t_k}{w}}^{\frac{t_{k+1}}{w}}[w^{\beta}|u-\ln x|^{\beta}+1]\;\widetilde\omega\left(f, w^{-\beta}\right)du\right)\\
&= \sum_{|w\ln x-t_k|\leq w/2}L(e^{-t_k}x^w)\,\psi\left(\left[1+\frac{w}{\Delta_{k}}\int_{\frac{t_k}{w}}^{\frac{t_{k+1}}{w}}w^{\beta}|u-\ln x|^{\beta}du\right]\;\widetilde\omega\left(f, w^{-\beta}\right)\right).
\end{split}
\end{equation*}
Since $\psi$ is concave and then subadditive, by (\ref{disPSI}) we have
\begin{equation*}
\begin{split}
I_{1,1}&\leq \sum_{|w\ln x-t_k|\leq w/2}L(e^{-t_k}x^w)\left[1+\frac{w}{\Delta_{k}}\int_{\frac{t_k}{w}}^{\frac{t_{k+1}}{w}}w^{\beta}|u-\ln x|^{\beta}du\right]\psi\left(\widetilde\omega\left(f, w^{-\beta}\right)\right)\\
&\leq M_{0,\Pi}(L)\,\psi\left(\widetilde\omega\left(f, w^{-\beta}\right)\right)\\
&+\psi\left(\widetilde\omega\left(f, w^{-\beta}\right)\right)\sum_{|w\ln x-t_k|\leq w/2}L(e^{-t_k}x^w)\frac{w}{\Delta_{k}}\int_{\frac{t_k}{w}}^{\frac{t_{k+1}}{w}}w^{\beta}|u-\ln x|^{\beta}du,
\end{split}
\end{equation*}
for $w>0$ sufficiently large. By using again (\ref{diseq}) and by exploiting the subadditivity of the function $|\cdot|^{\beta}$, with $0<\beta<1$, we can write
\begin{equation*}
\begin{split}
I_{1,1}&\leq \psi\left(\widetilde\omega\left(f, w^{-\beta}\right)\right)\left[M_{0,\Pi}(L)+\sum_{|w\ln x-t_k|\leq w/2} L(e^{-t_k}x^w)\left(w^{\beta}\left|\frac{t_{k}}{w}-\ln x\right|^{\beta}+\Delta^{\beta}\right)\right]\\
&\leq \psi\left(\widetilde\omega\left(f, w^{-\beta}\right)\right)\Biggl[M_{0,\Pi}(L)+\sum_{|w\ln x-t_k|\leq w/2} L(e^{-t_k}x^w)\left|w\ln x-t_{k}\right|^{\beta}\\
&\qquad\qquad\qquad\qquad\quad+\,\Delta^{\beta}\sum_{|w\ln x-t_k|\leq w/2} L(e^{-t_k}x^w)\Biggr]\\
&\leq \psi\left(\widetilde\omega\left(f, w^{-\beta}\right)\right)\left[M_{0,\Pi}(L)+M_{\beta,\Pi}(L)+\,\Delta^{\beta}M_{0,\Pi}(L)\right].
\end{split}
\end{equation*}
For what concerns $I_{1,2}$, we have
\begin{equation*}
\begin{split}
I_{1,2}&\leq \psi(2\left\|f\right\|_{\infty})\sum_{|w\ln x-t_k|> w/2}L(e^{-t_k}x^w)\\
&\leq \psi(2\left\|f\right\|_{\infty})\sum_{|w\ln x-t_k|> w/2}\frac{|w\ln x-t_k|^{\beta}}{|w\ln x-t_k|^{\beta}}L(e^{-t_k}x^w)\\
&\leq \left(\frac{2}{w}\right)^{\beta}\,\psi(2\left\|f\right\|_{\infty})\sum_{|w\ln x-t_k|> w/2}|w\ln x-t_k|^{\beta}L(e^{-t_k}x^w)\\
&\leq 2^{\beta+1}\;\psi(\left\|f\right\|_{\infty})\;w^{-\beta}\;M_{\beta,\Pi}(L).
\end{split}
\end{equation*}
\end{description}

\noindent 
Finally, we estimate $I_2$ (refer always to (\ref{split_cont})). Setting again $A_{j}:=\{x\in\R^+\,:\,0\leq |f(x)|<1/j\}$, we can rewrite $I_2$ as follows  
\begin{equation*}
\begin{split}
		I_{2}&=\left|\sum_{k\in\Z}\chi(e^{-t_k}x^w,f(x))-f(x)\right|\\
		&\leq \left|\sum_{k\in\Z}\chi(e^{-t_k}x^w,f(x)\mathbf{1}_{A_j}(x))-f(x)\mathbf{1}_{A_j}(x)\right|\\
		&+\left|\sum_{k\in\Z}\chi(e^{-t_k}x^w,f(x)\mathbf{1}_{\R^+\setminus A_j}(x))-f(x)\mathbf{1}_{\R^+\setminus A_j}(x)\right|.
\end{split}
\end{equation*}
Therefore, by condition ($\chi$4), there exist constants $M_1,M_2,\alpha>0$ such that
\begin{equation*}
\begin{split}
		I_2&\leq \mathcal{S}_{w}^{j}(x)\mathbf{1}_{A_j}(x)+|f(x)|\left|\frac{1}{|f(x)|}\sum_{k\in\Z}\chi(e^{-t_k}x^w,f(x)\mathbf{1}_{\R^+\setminus A_j}(x))-\mathbf{1}_{\R^+\setminus A_j}(x)\right|\\
		&\leq M_1 w^{-\alpha}+|f(x)|\mathcal{T}_{w}^{j}(x)\mathbf{1}_{\R^+\setminus A_j}(x)\\
		&\leq M_1 w^{-\alpha}+|f(x)|M_2 w^{-\alpha}\\
		&\leq M_1 w^{-\alpha}+M_2\left\|f\right\|_{\infty}w^{-\alpha},
	\end{split}
\end{equation*}
uniformly with respect to $x\in\R^+$, for sufficiently large $w>0$. Thus, the theorem is proved.
\end{proof}
Recalling that a function $f : \mathbb{R}^{+} \rightarrow \mathbb{R}$ is called \textit{locally log-H\"olderian} of order $0 < \nu \leq 1$ if
\begin{equation*}
\widetilde\omega(f, \delta) = \mathcal{O}(\delta^\nu), \quad \text{as } \delta \rightarrow 0^+,
\end{equation*}
we can state the following corollary.

\begin{corollary}
Let $f : \mathbb{R}^{+} \rightarrow \mathbb{R}$ be a locally log-H\"olderian function of order, $0<\nu\leq 1$. Suppose that $\psi$ of condition ($\chi$3) is concave and it satisfies the following assumption
\begin{equation*}
	\psi(u)=\mathcal{O}(u^q), 
\end{equation*}
as $u\to 0^{+}$, for some $0<q\leq 1$. Then, there exist constants $C_1, C_2>0$ such that
\begin{enumerate}
	\item if $L$ satisfies condition (L2) with $\beta\geq 1$:
\begin{equation*}
\left\|K^\chi_w f-f\right\|_{\infty}\leq C_1 w^{-l},
\end{equation*}
\item if $L$ satisfies condition ($L$2) with $0<\beta< 1$:
\begin{equation*}
\begin{split}
\left\|K^\chi_w f-f\right\|_{\infty}&\leq C_2 w^{-m},
\end{split}
\end{equation*}
\end{enumerate}
for sufficiently large $w>0$, with $l:=\min\{\nu q, \alpha\}$ and $m:=\min\{\nu\beta q, \beta, \alpha\}$, where $\alpha>0$ is the constant of condition ($\chi 4$).
\end{corollary}

We now conclude this section by presenting a Voronovskaja-type theorem for the nonlinear exponential Kantorovich sampling operators $(K_w^{\chi}f)_{w>0}$ by means of Mellin derivatives $\theta f$. This mathematical tool, introduced by Butzer and Jansche in \cite{butzer1997direct}, is defined as
\begin{equation*}
\theta f(x) := xf'(x), \quad x \in \mathbb{R}^+,
\end{equation*}
provided the usual derivative $f'(x)$ exists.\\

For the linear case, precise approximation formulas have been obtained using assumptions on the local moments of the kernel (see, \cite[Theorem 2.]{angamuthu2020direct}). In the case of nonlinear operators, we provide a partial estimate of the approximation error in terms of limsup, involving only the first derivative of the function $f$, in the Mellin-sense (see,\cite{bardaro2007voronovskaya}).\\

To achieve this, we assume the following additional condition on the function $L$ from ($\chi$3):
\begin{description}
\item {($L$3)} there exists $0<r\leq 1$, such that
\begin{equation*}
\lim_{w\to+\infty}\sum_{|t_k-w\ln x|>\gamma w} L(e^{-t_k}x^w)|t_k-w\ln x|^r=0,
\end{equation*}
uniformly with respect to $x\in\R^{+}$.
\end{description}

\begin{theorem}
Let $0<r\leq 1$, with $r<\alpha$, where $\alpha$ is the constant from condition ($\chi 4$). Suppose that ($\chi$3) is satisfied with $\psi(u) = u^{r}$, $u \in\R^+$, and that the function $L$ satisfies condition ($L$2) with the parameter $\beta=r$ and also ($L$3) with $r$. Then, for every $f\in \mathcal{C}(\R^{+})$ such that $(\theta f)(x)$ exists at a point $x\in \R^{+}$, we have
\begin{equation*}
\limsup_{w\to +\infty} w^{r}|(K_{w}^{\chi}f)(x)-f(x)|\leq \frac{\Delta^r |(\theta f)(x)|^r}{2^r} M_{0,\Pi}(L)+ |(\theta f)(x)|^r M_{r,\Pi}(L).\end{equation*}
\end{theorem}
\begin{proof}
In view of (\ref{split_cont}) and using the fact that $\psi(u)=u^{r}$ with $0<r<\min\{\alpha,1\}$, we have
\begin{equation*}
\begin{split}
	w^{r}|(K_{w}^{\chi}f)(x)-f(x)|&\leq  w^{r}\sum_{k\in\Z} L( e^{-t_k}x^w)\left[\frac{w}{\Delta_{k}}\int_{\frac{t_k}{w}}^{\frac{t_{k+1}}{w}}|f(e^u)-f(x)|du\right]^{r}\\
&+w^{r}\left|\sum_{k\in\Z}\chi\left ( e^{-t_k}x^w,f(x)\right )-f(x)\right|:=J_1+J_2.
\end{split}
\end{equation*}
Considering assumption ($\chi$4) and the fact that $r < \alpha$, it follows that the term $J_{2}$ tends to zero as $w\to+\infty$ uniformly with respect to $x\in\R^+$. \\
Now, we evaluate the term $J_1$. Since $f$ is Mellin differentiable at point $x$, there exists a bounded function $h$ such that $\lim_{y\to 1}h(y)=0$ and
\begin{equation*}
f(e^u)=f(x)+(\theta f)(x)(u-\ln x)+h\left(\frac{e^u}{x}\right)(u-\ln x).
\end{equation*}
Thus, using the concavity of the function $(\cdot)^{r}$ over $\R^{+}$, we get
\begin{equation*}
\begin{split}
J_1&\leq w^{r}\sum_{k\in\Z} L( e^{-t_k}x^w)\left[\frac{w}{\Delta_{k}}\int_{\frac{t_k}{w}}^{\frac{t_{k+1}}{w}}|(\theta f)(x) (u-\ln x)|du\right]^{r}\\
&+w^{r}\sum_{k\in\Z} L( e^{-t_k}x^w)\left[\frac{w}{\Delta_{k}}\int_{\frac{t_k}{w}}^{\frac{t_{k+1}}{w}}\left|h\left(\frac{e^u}{x}\right)(u-\ln x)\right|du\right]^{r}:=J_{1,1}+J_{1,2}.
\end{split}
\end{equation*}
Taking into account that $(\cdot)^{r}$ is non decreasing and concave, and therefore subadditive,
\begin{equation*}
\begin{split}
J_{1,1}&=w^{r}\sum_{k\in\Z} L( e^{-t_k}x^w)\left[\frac{w}{\Delta_{k}}|(\theta f)(x)|\int_{\frac{t_k}{w}}^{\frac{t_{k+1}}{w}}|u-\ln x|du\right]^{r}\\
&\leq w^{r}\sum_{k\in\Z} L( e^{-t_k}x^w)\left[\frac{w}{\Delta_{k}}|(\theta f)(x)|\int_{\frac{t_k}{w}}^{\frac{t_{k+1}}{w}}\left(\left|u-\frac{t_k}{w}\right|+\left|\frac{t_k}{w}-\ln x\right|\right)du\right]^{r}\\
&\leq w^{r}\sum_{k\in\Z} L( e^{-t_k}x^w)\left[|(\theta f)(x)|\left(\frac{\Delta}{2w}+\left|\frac{t_k}{w}-\ln x\right|\right)\right]^{r}\\
&\leq \frac{\Delta^r |(\theta f)(x)|^r}{2^r} M_{0,\Pi}(L)+ |(\theta f)(x)|^r M_{r,\Pi}(L).
\end{split}
\end{equation*}
To evaluate $J_{1,2}$, let $\varepsilon>0$ be fixed, then there exists $\gamma>0$ such that $|h(y)|\leq\varepsilon$ whenever $|y-1|\leq\gamma$. Moreover, fix $w'$ such that $\frac{\Delta}{w}<\frac{\gamma}{2}$ for every $w>w'$. Then,
\begin{equation*}
\begin{split}
J_{1,2}&=w^{r}\left\{\sum_{|t_k-w\ln x|<\frac{w\gamma}{2}}+\sum_{|t_k-w\ln x|\geq\frac{w\gamma}{2}}\right\} L( e^{-t_k}x^w)\left[\frac{w}{\Delta_{k}}\int_{\frac{t_k}{w}}^{\frac{t_{k+1}}{w}}\left|h\left(\frac{e^u}{x}\right)(u-\ln x)\right|du\right]^{r}\\
&=J_{1,2}'+J_{1,2}''
\end{split}
\end{equation*}
Let us consider $J_{1,2}'$. Since $u$ lies within $ \left[\frac{t_k}{w},\frac{t_{k+1}}{w}\right]$, in light of (\ref{disGamma}), and recalling the local equivalence between the continuity and the log-continuity, there holds
\begin{equation*}
	\left|\ln\left(\frac{e^u}{x}\right)\right|=|u-\ln x|\leq\gamma,
\end{equation*}
whenever $\left|\frac{t_k}{w}-\ln x\right|<\frac{\gamma}{2}$. Therefore, using the sub-additivity of $(\cdot)^r$ and Lemma \ref{lemma_preliminare} (ii), we readily derive that 
\begin{equation*}
\begin{split}
J_{1,2}'&\leq w^r \varepsilon^r \sum_{|t_k-w\ln x|<\frac{w\gamma}{2}}L(e^{-t_k}x^w)\left(\frac{\Delta}{2w}+\frac{\left|t_k-w\ln x\right|}{w}\right)^r\\
&\leq \left(\frac{\varepsilon\Delta}{2}\right)^r M_{0,\Pi}(L)+\varepsilon^r M_{r,\Pi}(L).	
\end{split}
\end{equation*}
Finally, using the boundedness of $h$, the sub-additivity of $(\cdot)^r$, Lemma \ref{lemma_preliminare} (ii)  and assumption ($L3$), we get that 
\begin{equation*}
\begin{split}
J_{1,2}''&\leq w^r \left\|h\right\|^{r}_{\infty} \sum_{|t_k-w\ln x|\geq\frac{w\gamma}{2}}L(e^{-t_k}x^w)\left(\frac{\Delta}{2w}+\frac{\left|t_k-w\ln x\right|}{w}\right)^r\\
&\leq  \left\|h\right\|^{r}_{\infty}\frac{\Delta^r}{2^r}\sum_{|t_k-w\ln x|\geq\frac{w\gamma}{2}}L(e^{-t_k}x^w) \;+ \;\left\|h\right\|^{r}_{\infty}\sum_{|t_k-w\ln x|\geq\frac{w\gamma}{2}}L(e^{-t_k}x^w)\left|t_k-w\ln x\right|^r\\
&\leq \varepsilon\left\|h\right\|^{r}_{\infty}\left(\frac{\Delta^r}{2^r}+1\right) .	
\end{split}	
\end{equation*}
This completes the proof.
\end{proof}

\section{Results in Mellin-Orlicz spaces $L^{\varphi}_{\mu}(\R^{+})$}\label{SectMellO}
It would be interesting to formulate the approximation results of Section \ref{MainSec} not only for log-uniformly continuous functions, but also for functions belonging to some general function spaces.\\

Let $\varphi:\mathbb{R}_{0}^{+}\to\mathbb{R}^{+}_{0}$ be a function. It is well-known that $\varphi$ is said to be a \textit{$\varphi$-function} if it satisfies the following conditions:
\begin{description}
\item ($\varphi$1) $\varphi$ is a non decreasing and continuous function;
\item ($\varphi$2) $\varphi(0)=0$, $\varphi(u)>0$ for every $u>0$;
\item ($\varphi$3) $\varphi(u) \rightarrow + \infty$ if $u \rightarrow + \infty$.
\end{description}
Now, in order to give the definition of the Mellin-Orlicz spaces, we introduce the notion of the modular functional $I^{\varphi}_{\mu}$ associated to the $\varphi$-function $\varphi$, defined by
\begin{equation*}
I^{\varphi }_{\mu}[f]:=\int_{0}^{+\infty}\varphi (\left | f(x) \right |)\frac{dx}{x},
\end{equation*}
for every $f\in \mathcal{M}(\mathbb{R}^{+})$. The Mellin-Orlicz space generated by $\varphi$ can now be defined by 
\begin{center}
$L^{\varphi}_{\mu}(\mathbb{R}^{+}):=$ \{$f\in \mathcal{M}(\mathbb{R}^{+}) \:\: : I^{\varphi}_{\mu}[\lambda f]< + \infty$, for some $\lambda>0$ \}.	
\end{center}
In Mellin-Orlicz spaces, different notions of convergence can be introduced. In this paper, we recall the most natural notion of convergence in this setting, that is called \textit{modular convergence}. A net of functions $(f_w)_{w>0}\subset L^{\varphi}_{\mu}(\mathbb{R}^{+})$ is modularly convergent to $f\in L^{\varphi}_{\mu}(\mathbb{R}^{+})$, if there exists $\lambda >0$ such that
\begin{equation*}
I^{\varphi}_{\mu}\left[\lambda (f_w-f) \right]	=\int_{0}^{+\infty}\varphi (\lambda \left | f_w(x)-f(x) \right |)\frac{dx}{x}\to 0,
\end{equation*}
as $w\to +\infty$.\\

To establish a modular convergence theorem in Mellin-Orlicz spaces, we begin by examining the modular convergence in $C_{c}(\mathbb{R}^+)$. The following lemma plays a crucial role in demonstrating this convergence in $C_{c}(\mathbb{R}^+)$.

\begin{lemma}\label{new_lemma}
Let $L$ be a function satisfying conditions ($L$1) and ($L$2). We have that, for every $\gamma>0$ and $\varepsilon>0$, there exists a sufficiently large $M>0$ such that
\begin{equation*}
\int_{|\ln x|>M} w L(e^{-t_k}x^w)\frac{dx}{x}<\varepsilon,
\end{equation*}
for sufficiently large $w>0$ and all the elements $|t_k|\leq \gamma w$.
\end{lemma}
\begin{proof}
Let $\varepsilon>0$ be fixed. Since $L\in L^{1}_{\mu}(\R^+)$, there exists a constant $\widetilde M$ such that
\begin{equation*}
	\int_{|\ln x|>\widetilde M} L(x)\frac{dx}{x}<\varepsilon.
\end{equation*}
Now, let $\gamma>0$, and for every $w\geq 1$, consider all integers $k$ such that $t_k/w\in [-\gamma,\gamma]$. Choose $M>0$ large enough such that $M-\gamma>\widetilde{M}$. By the change of variable $w\ln x-t_{k}=\ln y$, we obtain
\begin{equation*}
	\int_{|\ln x|>M} wL(e^{-t_k}x^w)\frac{dx}{x}\leq \int_{|\ln y|>(M-\gamma )w} L(y)\frac{dy}{y} \leq \int_{|\ln y|>\widetilde M}L(y)\frac{dy}{y}<\varepsilon,
\end{equation*}
which completes the proof of the lemma.
\end{proof}

Thus, we can state the following theorem.

\begin{theorem}\label{convCc}
Let $\varphi$ be a convex $\varphi$-function. For every $f\in C_c(\R^{+})$ and $\lambda>0$, there holds
\begin{equation*}
	\lim_{w\to +\infty} I^{\varphi}_{\mu}\left[\lambda\left(K_{w}^{\chi}f-f\right)\right]=0.
\end{equation*}
\end{theorem}
\begin{proof}
We have to prove that
\begin{equation*}
\lim_{w\to+\infty} I^{\varphi}_{\mu}\left[\lambda \left( K^{\chi}_{w}f-f \right )\right]=\lim_{w\to+\infty}\int_{0}^{+\infty} \varphi \left ( \lambda \left | (K^\chi_{w}f)(x)-f(x) \right | \right )\frac{dx}{x}=0,
\end{equation*}
for every $\lambda>0$, which is equivalent to show that the family $\left(\varphi \left ( \lambda \left | K^\chi_{w}f-f \right | \right )\right)_{w>0}$ converges to zero in $L^{1}_{\mu}(\mathbb{R}^n)$, for every $\lambda>0$. In order to establish this, we will exploit the Vitali convergence theorem, for $p=1$. Let now $\lambda>0$ be fixed. Let $f\in C_c(\R^+)$. By Theorem \ref{PointUnifTH} and the continuity of $\varphi$, it is easy to see that $\displaystyle\lim_{w\to +\infty}\varphi(\lambda\left\|K^\chi _w f-f\right\|_{\infty})=0$, for every $\lambda>0$. Then, for every fixed $\varepsilon>0$ there exists $\overline{w}>0$ such that for every $w\geq\overline{w}$, we have
\begin{equation*}
\varphi\left( \lambda\left|(K^\chi_{w}f)(x)-f(x)\right| \right)\leq \varphi \left(\lambda \left\| K^\chi_{w}f-f \right \|_{\infty}\right) < \varepsilon,	
\end{equation*}
for every $x\in\mathbb{R}^+$. Therefore,
\begin{equation*}
\mu \left ( \left \{ x\in\mathbb{R}^+\: :\: \varphi \left ( \lambda \left | (K^\chi_{w}f)(x)-f(x) \right | \right )>\varepsilon  \right \} \right )=\mu \left ( \emptyset  \right )=0,	
\end{equation*}
for every $w\geq\overline{w}$, where $\mu$ denotes again the logarithmic measure. It follows that $\left(\varphi\left(\lambda \left| K^{\chi}_{w}f-f \right|\right)\right)_{w>0}$ converges in measure to zero. \newline
Let now $\varepsilon>0$ be fixed and let $r>1$ be such that $\supp f\subset \left[\frac{1}{r},r\right]$. Then, if $t_{k}\notin [-w	\gamma,w\gamma]$, with $\gamma> \ln r +\Delta$, we have, for sufficiently large $w>0$, $\left[\frac{t_k}{w}, \frac{t_{k+1}}{w}\right]\cap \left[\frac{1}{r},r\right]=\emptyset$, and so
\begin{equation*}
	\int_{\frac{t_k}{w}}^{\frac{t_{k+1}}{w}} f(e^u)du=0.
\end{equation*}
By Lemma \ref{new_lemma}, there exists a constant $M>0$ (we can assume $M>\ln r$ without any loss of generality), such that 
\begin{equation*}
\int_{|\ln x|>M} w L(e^{-t_k}x^w)\frac{dx}{x}<\varepsilon,
\end{equation*}
for sufficiently large $w>0$ and $\left|t_{k}\right|\leq w\gamma$. Therefore, employing Jensen inequality (see, e.g., \cite{costarelli2015sharp}), taking into account that the measure $\frac{1}{M_{0,\Pi}(L)}\sum_{\left|t_{k}\right|\leq w\gamma}L(e^{-t_k}x^w)\leq 1$, and the Fubini-Tonelli theorem, we can conclude that
\begin{equation*}
\begin{split}
\int_{|\ln x|>M}&\varphi\left(\lambda\left|(K_w^\chi f)(x)-f(x)\right|\right)\frac{dx}{x}=\int_{|\ln x|>M} \varphi(\lambda|(K^{\chi}_w f)(x)|) \frac{dx}{x}\\
&\leq \int_{|\ln x|>M}\varphi\left(\lambda\sum_{\left|t_{k}\right|\leq w\gamma}\left|\chi\left ( e^{-t_k}x^w,\frac{w}{\Delta_{k}}\int_{\frac{t_k}{w}}^{\frac{t_{k+1}}{w}}f(e^u)du \right )\right|\right)\frac{dx}{x}\\
&\leq \int_{|\ln x|>M}\varphi\left(\lambda\sum_{\left|t_{k}\right|\leq w\gamma} L(e^{-t_k}x^w)\psi\left (\frac{w}{\Delta_{k}}\int_{\frac{t_k}{w}}^{\frac{t_{k+1}}{w}}|f(e^u)|du \right )\right)\frac{dx}{x}\\
&\leq \int_{|\ln x|>M}\varphi\left(\lambda\sum_{\left|t_{k}\right|\leq w\gamma} L(e^{-t_k}x^w)\psi\left (\left\|f\right\|_{\infty}\right)\right)\frac{dx}{x}\\
&\leq \sum_{\left|t_{k}\right|\leq w\gamma}\frac{\varphi(\lambda M_{0,\Pi}(L)\psi\left (\left\|f\right\|_{\infty}\right))}{w M_{0,\Pi}(L)}  \int_{|\ln x|>M}w L(e^{-t_k}x^w)\frac{dx}{x}\\
&<\frac{\varepsilon}{w M_{0,\Pi}(L)} \sum_{\left|t_{k}\right|\leq w\gamma}\varphi(\lambda M_{0,\Pi}(L)\psi\left (\left\|f\right\|_{\infty}\right)) .
\end{split}
\end{equation*}
The number of terms in the sum above corresponds to the count of $t_{k}/w$ within the interval $[-\gamma,\gamma]$, which does not exceed $2\left[\frac{w\gamma}{\delta}\right]+2$, where $\left[\cdot\right]$ represents the integer part. Thus,
\begin{equation*}
\int_{|\ln x|>M} \varphi(\lambda|(K^{\chi}_w f)(x)|) \frac{dx}{x}<\frac{2\varepsilon}{M_{0,\Pi}(L)} \left(\left[\frac{\gamma}{\delta}\right]+1\right)\varphi(\lambda m_{0,\Pi^n}(L)\psi\left (\left\|f\right\|_{\infty}\right))=:\varepsilon\cdot C,
\end{equation*}
for every $w\geq 1$. Therefore, for $\varepsilon>0$ there exists a set $E_{\varepsilon}=[-M, M]$ such that for every measurable set $F$, with $F\cap E_{\varepsilon}=\emptyset$, we have
\begin{equation*}
\begin{split}
	\int_{F} \varphi(\lambda|(K^\chi_w f)(x)-f(x)|) \frac{dx}{x}&=\int_{F} \varphi(\lambda|(K^\chi_w f)(x)|) \frac{dx}{x}\\
	&\leq \int_{|\ln x|>M} \varphi(\lambda|(K^\chi_w f)(x)|) \frac{dx}{x}< \varepsilon\cdot C .
\end{split}
\end{equation*}
Finally, let $B\subset\R^+$ be a measurable set with $\mu(B)<\varepsilon/s$, where
\begin{equation*}
s:=\max\{\varphi(2\lambda M_{0,\Pi}(L)\psi(\left\|f\right\|_{\infty})),\varphi(2\lambda\left\|f\right\|_{\infty})\},
\end{equation*}
$\left\|f\right\|_{\infty}\neq 0$. Using Remark \ref{remark_bounded} and the convexity of $\varphi$, in correspondence to $\varepsilon>0$ and for every $w>0$,
\begin{equation*}
\begin{split}
\int_{B} \varphi(\lambda|(K^\chi_w f)(x)&-f(x)|) \frac{dx}{x}\leq \frac{1}{2}\int_{B}\varphi(2\lambda|(K^\chi_w f)(x)|) \frac{dx}{x}+\frac{1}{2}\int_{B}\varphi(2\lambda| f(x)|) \frac{dx}{x}\\
&\leq \frac{1}{2}\int_{B}\varphi(2\lambda M_{0,\Pi}(L)\psi(\left\|f\right\|_{\infty})) \frac{dx}{x}+\frac{1}{2}\int_{B}\varphi(2\lambda\left\|f\right\|_{\infty}) \frac{dx}{x}\\
&\leq \int_{B}s\, \frac{dx}{x}=\mu(B)\,s<\varepsilon.
\end{split}
\end{equation*}
It follows that the integrals $\displaystyle\int_{(\cdot)} \varphi(\lambda|(K^\chi_w f)(x)-f(x)|) \frac{dx}{x}$ are equi-absolutely continuous, and the proof is now complete.
\end{proof}

In order to show that the nonlinear exponential Kantorovich sampling operators are well-defined in the Mellin-Orlicz spaces, we introduce the following growth condition (H), which is commonly used for approximation by means of nonlinear operators (see \cite{MonografiaBMV}).\\

Let $\varphi$ be a fixed $\varphi$-function. We assume the existence of a $\varphi$-function $\eta$ such that, for every $\lambda \in (0,1)$, there exists a constant $C_\lambda \in (0,1)$ satisfying
\begin{equation}
\varphi(C_\lambda \psi (u)) \leq \eta(\lambda u),
\tag{H}
\end{equation}
for every $u\in\R^{+}_{0}$, where $\psi$ is the function associated with condition ($\chi$3).\\

Under this assumption, we can prove the following modular inequality.

\begin{theorem}\label{stimamModOrlicz}
Let $\varphi$ be a convex $\varphi$-function satisfying condition (H) with $\eta$ convex. Then, for any $f,g \in L^{\eta}_{\mu}(\mathbb{R}^+)$, there exist $\lambda\in (0,1)$ and a constant $0<c\leq C_\lambda/M_{0,\Pi}(L)$ such that
\begin{equation*}
I^{\varphi}_{\mu}[c (K^\chi_w f- K^\chi_w g)]\leq \frac{\left\|L \right\|_{1,\mu}}{\delta M_{0,\Pi}(L)}I^{\eta}_{\mu}[\lambda (f-g)].
\end{equation*}
\end{theorem}
\begin{proof}
Let $f,g \in L^{\eta}_{\mu}(\mathbb{R}^+)$. For $x\in\R^+$, by applying condition ($\chi$3), we have
\begin{equation*}
\begin{split}
|(K^\chi_w f)(x)- (K^\chi_w g)(x)|\leq \sum_{k\in\Z}L(e^{-t_k}x^w)\psi\left(\frac{w}{\Delta_{k}}\int_{\frac{t_k}{w}}^{\frac{t_{k+1}}{w}}|f(e^u)-g(e^u)|du\right ).
\end{split}
\end{equation*}
Since $f-g\in L^{\eta}_\mu(\R^n)$, there exists $\lambda>0$ (that can be considered $\lambda\in (0,1)$ without any loss of generality) such that $I^{\eta}_{\mu}[\lambda(f-g)]<+\infty$. Then, we can choose $c>0$ such that $c \,M_{0,\Pi}(L)\leq C_\lambda$, where $C_{\lambda}\in (0,1)$ is the parameter arising from condition (H).\newline
Therefore, applying Jensen inequality twice and Fubini-Tonelli theorem, together with the changes of variable $w\ln x-t_k=\ln y$ and $u=\ln t$, we get
\begin{equation*}
\begin{split}
I^{\varphi}_{\mu}&[c (K^\chi_w f- K^\chi_w g)]=\int_{0}^{+\infty}\varphi(c |(K^\chi_w f)(x)- (K^\chi_w g)(x)|)\frac{dx}{x}\\
&\leq \int_{0}^{+\infty}\varphi\left(c\sum_{k\in\Z}L(e^{-t_k}x^w)\psi\left(\frac{w}{\Delta_{k}}\int_{\frac{t_k}{w}}^{\frac{t_{k+1}}{w}}|f(e^u)-g(e^u)|du \right )\right)\frac{dx}{x}\\
&\leq \frac{1}{M_{0,\Pi}(L)}\sum_{k\in\Z}\varphi\left(c\,M_{0,\Pi}(L)\psi\left(\frac{w}{\Delta_{k}}\int_{\frac{t_k}{w}}^{\frac{t_{k+1}}{w}}|f(e^u)-g(e^u)|du \right )\right)\int_{0}^{+\infty}L( e^{-t_k}x^w)\frac{dx}{x}
\end{split}
\end{equation*}
\begin{equation*}
\begin{split}
&\leq \frac{\left\|L\right\|_{1,\mu}}{w\,M_{0,\Pi}(L)}\sum_{k\in\Z}\varphi\left(C_\lambda\psi\left(\frac{w}{\Delta_{k}}\int_{\frac{t_k}{w}}^{\frac{t_{k+1}}{w}}|f(e^u)-g(e^u)|du \right )\right)\\
&\leq \frac{\left\|L\right\|_{1,\mu}}{w\,M_{0,\Pi}(L)}\sum_{k\in\Z}\eta\left(\lambda\frac{w}{\Delta_{k}}\int_{\frac{t_k}{w}}^{\frac{t_{k+1}}{w}}|f(e^u)-g(e^u)|du \right)\\
&\leq \frac{\left\|L\right\|_{1,\mu}}{w\,M_{0,\Pi}(L)}\sum_{k\in\Z}\frac{w}{\Delta_{k}}\int_{\frac{t_k}{w}}^{\frac{t_{k+1}}{w}}\eta\left(\lambda|f(e^u)-g(e^u)| \right)du\\
&=\frac{\left\|L\right\|_{1,\mu}}{w\,M_{0,\Pi}(L)}\sum_{k\in\Z}\frac{w}{\Delta_{k}}\int_{e^\frac{t_k}{w}}^{e^\frac{t_{k+1}}{w}}\eta\left(\lambda|f(t)-g(t)| \right)\frac{dt}{t}\\
&\leq \frac{\left\|L\right\|_{1,\mu}}{\delta\,M_{0,\Pi}(L)}I^{\eta}_{\mu}[\lambda(f-g)].
\end{split}
\end{equation*}
The proof is now complete.
\end{proof}

We finally prove the main theorem of this section.
\begin{theorem}\label{convModOrlicz}
Let $\varphi$ be a convex $\varphi$-function satisfying condition (H) with $\eta$ convex. If $f\in L^{\varphi+\eta}_{\mu}(\mathbb{R}^+)$, then there exists $c>0$ such that
\begin{equation*}
\lim_{w\to +\infty } I^{\varphi }_{\mu}[c(K_{w}f-f)]=0.	
\end{equation*}
\end{theorem}
\begin{proof}
Let $f\in L^{\varphi+\eta}_{\mu}(\mathbb{R}^+)$. By the density theorem (see Theorem 9.1 of \cite{MonografiaBMV}), there exists a $\lambda\in (0,1)$ such that, for every $\varepsilon >0$ there exists a function $g\in C_c(\R^+)$ such that $I^{\varphi+\eta}_{\mu}[\lambda(f-g)]<\varepsilon$. Now, we can fix a constant $c>0$ such that 
\begin{equation*}
	c\leq\min \left\{\frac{C_\lambda}{3M_{0,\Pi}(L)},\frac{\lambda}{3}\right\},
\end{equation*}
where $C_\lambda$ is the constant of condition (H). By the properties of the modular $I^{\varphi}_{\mu}$ and Theorem \ref{stimamModOrlicz}, we can write
\begin{equation*}
\begin{split}
I^{\varphi }_{\mu}[c(K_{w}f-f)]&\leq I^{\varphi }_{\mu}[3c(K_{w}f-K_w g)]+I^{\varphi }_{\mu}[3c(K_{w}g-g)]+I^{\varphi }_{\mu}[3c(f-g)]\\
&\leq \frac{\left\|L \right\|_{1,\mu}}{\delta M_{0,\Pi}(L)}I^{\eta}_{\mu}[\lambda (f-g)]+I^{\varphi }_{\mu}[\lambda(K_{w}g-g)]+I^{\varphi }_{\mu}[\lambda(f-g)].
\end{split}
\end{equation*}
Let $\kappa:=\max\left\{\frac{\left\|L \right\|_{1,\mu}}{\delta M_{0,\Pi}(L)},1\right\}$, we have
\begin{equation*}
\begin{split}
	I^{\varphi }_{\mu}[c(K_{w}f-f)]&\leq\kappa I^{\varphi +\eta}_{\mu}[\lambda(f-g)]+I^{\varphi }_{\mu}[\lambda(K_{w}g-g)]\\
	&\leq \kappa\varepsilon + I^{\varphi }_{\mu}[\lambda(K_{w}g-g)].
\end{split}
\end{equation*}
Hence, the assertions follows from Theorem \ref{convCc}.
\end{proof}
As the Mellin-Lebesgue spaces are particular instances of the Mellin-Orlicz spaces, with a $\varphi$-function of power type ($\varphi(u)=u^p$, with $1\leq p<+\infty$), Theorem \ref{convModOrlicz} also provides a convergence result within these contexts.\\

\section{Quantitative estimates in Mellin-Orlicz spaces $L^{\varphi}_{\mu}(\R^{+})$}
Here, we aim to obtain quantitative estimates in Mellin-Orlicz spaces. To achieve this, we need to introduce a suitable modulus of smoothness within such spaces. This concept is based on the Mellin translation operator $t_h f$, where $h\in\R^+$ and $(t_h f)(x):=f(hx)$, $x\in\R^+$. Consequently, for any fixed $f\in L^{\varphi}_{\mu}(\mathbb{R}^{+})$, we define the \textit{log-modulus of smoothness} in the Mellin-Orlicz spaces $L^{\varphi}_{\mu}(\mathbb{R}^{+})$, with respect to the modular $I^{\varphi}_{\mu}$, as follows
\begin{equation*}
\widetilde\omega_{\varphi,\mu}(f,\delta):=\sup_{\left | \ln t \right |\leq\delta}	I^{\varphi}_{\mu}\left[f(\cdot\, t)-f(\cdot)\right]=\sup_{\left | \ln t \right |\leq \delta}\int_{0}^{+\infty}\varphi(|f(st)-f(s)|)\frac{ds}{s},
\end{equation*}
with $\delta>0$. Moreover, it is well-known that for every $f\in L^{\varphi}_{\mu}(\mathbb{R}^{+})$ there exists $\lambda >0$ such that $\widetilde\omega_{\varphi,\mu}(\lambda f, \delta)\to 0$, as $\delta \to 0^{+}$ (see, \cite{MonografiaBMV}).\newline\newline

To achieve our goals, we need to require a different version of condition ($\chi$4). This condition (inspired from the one given in \cite{CV2015bis}) is stated as follows
\begin{description}
\item ($\chi$4*) there exists $\alpha>0$, such that
\begin{equation*}
\mathcal{T}_w(x):=\sup_{u\neq 0} \left |\frac{1}{u} \sum _{k\in\mathbb{Z}}\chi \left ( e^{-t_k}x^w,u \right )-1 \right |= \mathcal{O}\left (w^{-\alpha}  \right ),
\end{equation*}
as $w\to +\infty$, uniformly with respect to $x\in\R^{+}$.
\end{description}
Therefore, from now on, we will consider nonlinear kernels that satisfy condition ($\chi$4*) in place of condition ($\chi$4).\newline

Hence, denoting by $\tau$ the characteristic function of the set $[1, e]$, we can state the following quantitative estimates as desired.
\begin{theorem}\label{mainstime}
Let $\Pi=(t_k)_{k\in\Z}$ be a sequence of real numbers such that $\Delta_k=t_{k+1}-t_k=1$, for every $k\in\Z$. Suppose that $\varphi$ is a convex $\varphi$-function which satisfies condition (H) with $\eta$ convex, $f\in L^{\varphi +\eta}_{\mu}(\mathbb{R}^+)$ and also for any fixed $0<\gamma <1$, we have
\begin{equation}\label{e3.1}
w \int_{|\ln y|>\frac{1}{w^\gamma }}L(y^w)\frac{dy}{y}\leq M_{3}w^{-\gamma_{0}},	
\end{equation}
as $w\to +\infty$, for suitable positive constants $M_3$, $\gamma_0$ depending on $\alpha$ and $L$. Then, there exist $\nu>0$, $\lambda >0$ and further parameter $\lambda_0 >0$ such that
\begin{equation*}\begin{split}
I^\varphi_\mu [\nu (K^\chi_w f-f)]&\leq \frac{\left \| L \right \|_{1,\mu}M_{0,\Pi}(\tau)}{3  M_{0, \Pi}(L)}\widetilde\omega_{\eta ,\mu} \left ( \lambda f,\frac{1}{w^\gamma } \right )+\frac{M_3 M_{0,\Pi}(\tau)I^{\eta }_{\mu}[\lambda _0 f]}{3  M_{0, \Pi}(L)}w^{-\gamma _0}\\
&+\frac{1}{3}	\widetilde\omega_{\eta ,\mu} \left (\lambda f, \frac{1}{w} \right )+\frac{I^{\varphi }_{\mu}[\lambda _0 f]}{3}w^{-\alpha},	
\end{split}\end{equation*}
for every sufficiently large $w>0$, where $M_{0,\Pi}(L)<+\infty$ by (i) of Lemma \ref{lemma_preliminare}, $M_{0,\Pi}(\tau)<+\infty$, since $\tau$ is bounded and with compact support, and $\alpha>0$ is the constant of condition ($\chi 4$*).
\end{theorem}
\begin{proof}
Let $\lambda_0$ such that $I^{\varphi}_{\mu}[\lambda_0 f]<+\infty$. Further, we also fix $\lambda >0$ such that
\begin{equation*}
\lambda < \min\left \{ 1, \frac{\lambda _0}{2} \right \}.	
\end{equation*}
In correspondence to $\lambda$, by condition (H), we know that there exists $C_{\lambda}\in (0,1)$ such that $\varphi( C_{\lambda} \psi(u))\leq \eta(\lambda u)$, $u\in\mathbb{R}^{+}_{0}$, while by ($\chi$4*), there exist constants $\alpha$, $M>0$ such that $\mathcal{T} _w (x)\leq M w^{-\alpha}$, uniformly with respect to $\underline{x}\in\mathbb{R}^n$, for sufficiently large $w>0$. Now, we choose $\nu >0$ such that
\begin{equation*}
\nu\leq\min\left\{\frac{C_{\lambda}}{3M_{0,\Pi}(L)},\frac{\lambda_0}{3M}\right\}.	
\end{equation*}
Taking into account that $\varphi$ is convex and non-decreasing, for $\nu>0$, we can write
\begin{equation*}\begin{split}
&I^{\varphi }_{\mu}\left [ \nu \left ( K^\chi_w f -f \right ) \right ]=\int_{0}^{+\infty}\varphi \left ( \nu \left | (K^\chi_w f)(x)-f(x) \right | \right )\frac{dx}{x}\\
&\leq \frac{1}{3}\Biggl\{\int_{0}^{+\infty}\varphi \left ( 3\nu \left | (K^\chi_w f)(x)-\sum_{k\in\mathbb{Z}}\chi\left ( e^{-t_k}x^w, w\int_{\frac{t_k}{w}}^{\frac{t_{k+1}}{w}}f\left ( e^{u-\frac{t_k}{w}}x \right ) du\right ) \right | \right )\frac{dx}{x}\\
&+\int_{0}^{+\infty}\varphi \left ( 3\nu \left | \sum_{k\in\mathbb{Z}}\chi\left ( e^{-t_k}x^w, w\int_{\frac{t_k}{w}}^{\frac{t_{k+1}}{w}}f\left ( e^{u-\frac{t_k}{w}}x \right ) du\right ) -\sum_{k\in\mathbb{Z}}\chi\left (e^{-t_k}x^w,f(x) \right)\right | \right )\frac{dx}{x}\\
&+\int_{0}^{+\infty}\varphi \left ( 3\nu \left | \sum_{k\in\mathbb{Z}}\chi\left ( e^{-t_k}x^w,f(x) \right)-f(x)\right | \right )\frac{dx}{x} \Biggr\}=:I_1 +I_2 +I_3.
\end{split}\end{equation*}
Now, we estimate $I_1$. Applying condition ($\chi$3), we have
\begin{equation*}\begin{split}
&3I_1 = \int_{0}^{+\infty}\varphi \left ( 3\nu \left | (K^\chi_w f)(x)-\sum_{k\in\mathbb{Z}}\chi\left ( e^{-t_k}x^w, w\int_{\frac{t_k}{w}}^{\frac{t_{k+1}}{w}}f\left ( e^{u-\frac{t_k}{w}}x \right ) du\right ) \right | \right )\frac{dx}{x}\\
&\leq \int_{0}^{+\infty}\varphi \left ( 3\nu \sum_{k\in\mathbb{Z}} L(e^{-t_k}x^w) \psi\left(\left | w\int_{\frac{t_k}{w}}^{\frac{t_{k+1}}{w}}f\left ( e^u\right )- f\left ( e^{u-\frac{t_k}{w}}x \right ) du\right|\right )  \right )\frac{dx}{x}.
\end{split}\end{equation*}
Using Jensen inequality twice, the changes of variable $v=e^u$ and $\ln y=\ln x-\frac{t_{k}}{w}$, the fact that $\Delta_k=1$ for every $k\in\Z$, condition (H) and Fubini-Tonelli theorem, we obtain
\begin{equation*}\begin{split}
3I_1 &\leq \frac{1}{M_{0,\Pi}(L)}\int_{0}^{+\infty}\sum_{k\in\mathbb{Z}}L( e^{-t_k}x^w) \varphi \left ( 3\nu M_{0,\Pi}(L) \psi\left(w\int_{\frac{t_k}{w}}^{\frac{t_{k+1}}{w}}\left | f\left ( e^u\right )- f\left ( e^{u-\frac{t_k}{w}}x \right )\right| du\right )  \right )\frac{dx}{x}\\
&\leq\frac{1}{M_{0,\Pi}(L)}\sum_{k\in\mathbb{Z}}\int_{0}^{+\infty}L( e^{-t_k}x^w)\varphi \left ( C_{\lambda} \psi\left(w\int_{\frac{t_k}{w}}^{\frac{t_{k+1}}{w}}\left | f\left ( e^u\right )- f\left ( e^{u-\frac{t_k}{w}}x \right )\right| du\right )  \right )\frac{dx}{x}\\
&\leq \frac{1}{M_{0,\Pi}(L)}\sum_{k\in\mathbb{Z}}\int_{0}^{+\infty}L( e^{-t_k}x^w)\eta \left ( \lambda w\int_{\frac{t_k}{w}}^{\frac{t_{k+1}}{w}}\left | f\left ( e^u\right )- f\left ( e^{u-\frac{t_k}{w}}x \right )\right| du  \right )\frac{dx}{x}\\
&\leq \frac{1}{M_{0,\Pi}(L)}\sum_{k\in\mathbb{Z}}\int_{0}^{+\infty}L( e^{-t_k}x^w)w\int_{\frac{t_k}{w}}^{\frac{t_{k+1}}{w}}\eta \left ( \lambda \left |f\left ( e^u\right )- f\left ( e^{u-\frac{t_k}{w}}x \right )\right|\right ) du \frac{dx}{x}
\end{split}\end{equation*}
\begin{equation*}\begin{split}
&= \frac{1}{M_{0,\Pi}(L)}\sum_{k\in\mathbb{Z}}\int_{0}^{+\infty}L( e^{-t_k}x^w)w\int_{e^\frac{t_k}{w}}^{e^\frac{t_{k+1}}{w}}\eta \left ( \lambda\left |f\left ( v\right )- f\left ( ve^{-\frac{t_k}{w}}x  \right )\right|\right )\frac{dv}{v}\frac{dx}{x}\\
&= \frac{1}{M_{0,\Pi}(L)}\sum_{k\in\mathbb{Z}}\int_{0}^{+\infty}L( y^w)w\int_{0}^{+\infty}\eta \left ( \lambda\left |f\left ( v\right )- f\left ( vy\right )\right|\right )\tau(e^{-t_k}v^w)\frac{dv}{v}\frac{dy}{y}\\
&= \frac{1}{ M_{0,\Pi}(L)}\int_{0}^{\infty} w L( y^w) \int_{0}^{\infty}\eta \left ( \lambda\left |f\left ( v\right )- f\left ( vy \right )\right|\right )\sum_{k\in\mathbb{Z}}\tau(e^{-t_k}v^w) \frac{dv}{v}\frac{dy}{y}\\
&\leq \frac{M_{0,\Pi}(\tau)}{ M_{0,\Pi}(L)}\int_{0}^{+\infty} w L(y^w) \int_{0}^{+\infty}\eta \left ( \lambda\left |f\left ( v\right )- f\left ( vy\right )\right|\right )\frac{dv}{v}\frac{dy}{y}\\
&=\frac{M_{0,\Pi}(\tau)}{ M_{0,\Pi}(L)}\int_{0}^{+\infty} w L( y^w) \;I^{\eta} _{\mu}\left[ \lambda \left( f(\cdot)- f(\cdot \, y) \right)\right ] \frac{dy}{y},
\end{split}\end{equation*}
where the constant $M_{0,\Pi}(\tau)<+\infty$, since $\tau$ is bounded and with compact support (see, e.g., \cite{CV2015bis, CCNV2022}). 
Now, let $0<\gamma<1$ be fixed. We now split the above integral as follows
\begin{equation*}
\begin{split}
	&\frac{w M_{0,\Pi}(\tau)}{ M_{0,\Pi}(L)}\Biggl\{\int_{|\ln y|\leq\frac{1}{w^{\gamma }}}+\int_{|\ln y|>\frac{1}{w^{\gamma }}}\Biggr\} L( y^w)\; I^{\eta} _{\mu}\left[ \lambda \left( f(\cdot)- f(\cdot \, y) \right)\right ] \frac{dy}{y}\\
	&\qquad\qquad\qquad\qquad\qquad\qquad\qquad=:I_{1,1}+I_{1,2}.
\end{split}
\end{equation*}
For $I_{1,1}$, one has
\begin{equation*}\begin{split}
I_{1,1}&\leq \frac{M_{0,\Pi}(\tau)}{ M_{0,\Pi}(L)}\;\widetilde\omega_{\eta,\mu}\left(\lambda f, \frac{1}{w^{\gamma}} \right) w\int_{|\ln y|\leq\frac{1}{w^{\gamma }}} L( y^w) \frac{dy}{y}\\
&\leq \frac{M_{0,\Pi}(\tau)}{ M_{0,\Pi}(L)}\;\widetilde\omega_{\eta,\mu}\left(\lambda f, \frac{1}{w^{\gamma}} \right) \left \| L \right \|_{1,\mu}.
\end{split}\end{equation*}
On the other hand, taking into account that $\eta$ is convex, for $I_{1,2}$ we can write
\begin{equation*}
I_{1,2}\leq \frac{w M_{0,\Pi}(\tau)}{ M_{0,\Pi}(L)}\int_{|\ln y|>\frac{1}{w^{\gamma }}} L( y^w) \frac{1}{2}\left(I^{\eta}_\mu[2\lambda f(\cdot)]+I^{\eta}_\mu[2\lambda f(\cdot\, y)]\right) \frac{dy}{y}.
\end{equation*}
Now, observing that
\begin{equation*}
I^{\eta}_\mu[2\lambda f(\cdot)]=I^{\eta}_\mu[2\lambda f(\cdot \, y)],
\end{equation*}
for every fixed $y$, using (\ref{e3.1}), we finally get
\begin{equation*}\begin{split}
I_{1,2}&\leq \frac{w M_{0,\Pi}(\tau)}{M_{0,\Pi}(L)}\int_{|\ln y|>\frac{1}{w^{\gamma }}} L( y^w)\; I^{\eta}_\mu[2\lambda f] \frac{dy}{y}\\
&\leq \frac{M_{0,\Pi}(\tau)}{ M_{0,\Pi}(L)}\; I^{\eta}_\mu[\lambda_0 f] \; M_3 w^{-\gamma_0},
\end{split}\end{equation*}
for $w > 0$ sufficiently large and for $M_3 >0$.\\
Now we can proceed estimating $I_2$. Using the assumption ($\chi$3) we immediately have
\begin{equation*}\begin{split}
3 I_2 &= \int_{0}^{+\infty}\varphi \left ( 3\nu \left | \sum_{k\in\mathbb{Z}}\chi\left ( e^{-t_k}x^w, w\int_{\frac{t_k}{w}}^{\frac{t_{k+1}}{w}}f\left ( e^{u-\frac{t_k}{w}}x \right ) du\right ) -\sum_{k\in\mathbb{Z}}\chi\left (e^{-t_k}x^w,f(x) \right)\right | \right )\frac{dx}{x}\\
&\leq \int_{0}^{+\infty}\varphi \left ( 3\nu \sum_{k\in\mathbb{Z}}L( e^{-t_k}x^w) \; \psi\left(w\int_{\frac{t_k}{w}}^{\frac{t_{k+1}}{w}}\left|f\left ( e^{u-\frac{t_k}{w}}x \right ) - f(x) \right|du\right ) \right )\frac{dx}{x}
\end{split}\end{equation*}
Now, by the change of variable $\ln{y}=u-\frac{t_{k}}{w}$, we have
\begin{equation*}
3I_2\leq \int_{0}^{+\infty}\varphi \left ( 3\nu \sum_{k\in\mathbb{Z}}L( e^{-t_k}x^w) \; \psi\left(w\int_{1}^{e^{1/w}}\left|f\left ( xy\right)- f(x)\right|\frac{dy}{y}\right ) \right )\frac{dx}{x}.
\end{equation*}
Hence, taking into account that $w\int_{1}^{e^{1/w}}\frac{dy}{y}=1$, applying Jensen inequality twice, recalling that $3\nu M_{0, \Pi}(L)\leq C_{\lambda}$, using condition (H) and Fubini-Tonelli theorem, we get
\begin{equation*}\begin{split}
3I_2&\leq \frac{1}{M_{0,\Pi}(L)} \int_{0}^{+\infty} \sum_{k\in\mathbb{Z}}L( e^{-t_k}x^w) \;\varphi \left ( 3\nu M_{0,\Pi}(L)\; \psi\left(w\int_{1}^{e^{1/w}}\left|f\left ( xy\right)- f(x)\right|\frac{dy}{y}\right ) \right )\frac{dx}{x}\\
&\leq \frac{1}{M_{0,\Pi}(L)} \int_{0}^{+\infty} \sum_{k\in\mathbb{Z}}L( e^{-t_k}x^w) \;\eta\left ( \lambda w\int_{1}^{e^{1/w}}\left|f\left ( xy\right)- f(x)\right|\frac{dy}{y} \right )\frac{dx}{x}\\
&\leq \frac{1}{M_{0,\Pi}(L)}\int_{0}^{+\infty} \sum_{k\in\mathbb{Z}}L( e^{-t_k}x^w) w\int_{1}^{e^{1/w}}\eta\left ( \lambda\left|f\left ( xy \right)- f(x)\right|\right )\frac{dy}{y}\frac{dx}{x}\\
&\leq w \int_{0}^{+\infty} \int_{1}^{e^{1/w}}\eta\left ( \lambda\left|f\left ( xy\right)- f(x)\right|\right )\frac{dy}{y}\frac{dx}{x}\\
&=w\int_{1}^{e^{1/w}}I^{\eta}_\mu\left [ \lambda\left(f\left ( \cdot \, y \right)- f(\cdot)\right)\right]\frac{dy}{y}\\
&\leq \widetilde\omega_{\eta,\mu}\left(\lambda f, \frac{1}{w}\right) w \int_{1}^{e^{1/w}}\frac{dy}{y}\\
&= \widetilde\omega_{\eta,\mu}\left(\lambda f, \frac{1}{w}\right).
\end{split}\end{equation*}
For $I_3$, denoted by $A_{0}\subseteq \mathbb{R}^{+}$ the set of all points of $\mathbb{R}^{+}$ for which $f\neq 0$ a.e., we can write, by condition ($\chi$2) and using the fact that $\varphi(0)=0$, what follows
\begin{equation*}\begin{split}
3I_3&=\int_{A_0}\varphi \left ( 3\nu \left | \sum_{k\in\mathbb{Z}}\chi\left ( e^{-t_k}x^w,f(x) \right)-f(x)\right | \right )\frac{dx}{x}\\
&=\int_{A_0}\varphi \left ( 3\nu |f(x)|\left |\frac{1}{f(x)} \sum_{k\in\mathbb{Z}}\chi\left ( e^{-t_k}x^w,f(x) \right)-1\right | \right )\frac{dx}{x}\\
&\leq \int_{A_0}\varphi \left ( 3\nu |f(x)|\mathcal{T}_{w}(x) \right )\frac{dx}{x}.
\end{split}\end{equation*}
By the convexity of $\varphi$ and condition ($\chi$4*), we have
\begin{equation*}\begin{split}
3I_{3}&\leq \int_{A_0}\varphi \left ( 3\nu M w^{-\alpha}|f(x)|\right )\frac{dx}{x}\\
&\leq  w^{-\alpha}\int_{\mathbb{R}^{+}}\varphi \left ( 3\nu M |f(x)|\right )\frac{dx}{x}\\
&\leq w^{-\alpha} \int_{\mathbb{R}^{+}}\varphi \left ( \lambda_0 |f(x)|\right )\frac{dx}{x}\\
&=w^{-\alpha} I^{\varphi}_{\mu}\left[  \lambda_{0} f\right ],
\end{split}\end{equation*}
for positive constants $M$ and $\alpha$. This completes the proof.
\end{proof}

Now we, recall the definition of the Lipschitz (log-H\"olderian) classes in Mellin-Orlicz spaces $L^{\varphi}_{\mu}(\mathbb{R}^+)$. We define by $Lip_{\varphi,\mu}(\nu)$, $0<\nu\leq 1$, as follows
\begin{equation*}
Lip_{\varphi,\mu} (\nu):=\{f\in L^{\varphi}_{\mu}(\mathbb{R}^+) \;:\; \exists \lambda >0 ,\text{  } \widetilde\omega_{\nu,\mu}(\lambda f,\delta)=\mathcal{O}(\delta^{\nu}), \text{ as } \delta\to 0\}
\end{equation*}
From Theorem \ref{mainstime}, we obtain the following corollary.
\begin{corollary}\label{c3.1}
Under the assumptions of Theorem \ref{mainstime} with $0<\gamma<1$ and for any $f\in Lip_{\eta,\mu}(\nu)$, $0<\nu\leq 1$, there exist $C>0$ and $c>0$ such that
\begin{equation*}
	I^{\varphi}_\mu\left[c \left(K^\chi_w f -f\right)\right]\leq C w^{-l},
\end{equation*}
for sufficiently large $w>0$, with $l:=\min\{\gamma \nu, \gamma_0, \alpha\}$, where $\alpha >0$ is the constant of condition ($\chi 4$).
\end{corollary}

Finally, to support the theory with practical examples in the nonlinear context, we define kernel functions as follows
\begin{equation*}
\chi (e^{-t_k}x^w,u) = L(e^{-t_k}x^w) \, g_{w}(u),	
\end{equation*}
 where $(g_w)_{w>0}$, $g_w: \mathbb{R}\to\mathbb{R}$ is a family of functions satisfying $g_{w}(u)\to u$ uniformly as $w\to+\infty$ and such that \begin{equation*}
\left|g_w (u) - g_w (v)\right|	\leq \psi(|u-v|),
\end{equation*}
for every $u,v\in\mathbb{R}$ and $w>0$. For concrete examples of function sequences $(g_w)_{w>0}$ and specific functions $L$, we refer the readers to \cite{chapter_costa, costarelli2024comparison}.

\section*{Acknowledgment}
{\small The authors are members of the Gruppo Nazionale per l'Analisi Matematica, la Probabilit\`{a} e le loro Applicazioni (GNAMPA) of the Istituto Nazionale di Alta Matematica (INdAM), of the gruppo UMI (Unione Matematica Italiana) T.A.A. (Teoria dell' Approssimazione e Applicazioni)}, and of the network RITA (Research ITalian network on Approximation).

\section*{Funding}
{\small The authors have been supported in this research by the project PRIN 2022: ``AI- and DIP-Enhanced DAta Augmentation for Remote Sensing of Soil Moisture and Forest Biomass (AIDA)" funded by the European Union under the Italian National Recovery and Resilience Plan (NRRP) of NextGenerationEU, under the Italian Ministry of University and Research (MUR) (Project Code: 20229FX3B9, CUP: J53D23000660001).}

\section*{Conflict of interest/Competing interests}
{\small The authors declare that they have no conflict of interest and competing interest.}

\section*{Availability of data and material and Code availability}
{\small Not applicable.}

%--------------REFERENCES------------------------------------------------

%-------------------------------------------------------------------------
\end{document}